\documentclass[11pt]{article}
\usepackage{epsfig}
\usepackage{t1enc}
    \usepackage[latin1]{inputenc}
    \usepackage[english]{babel}
        \usepackage{latexsym}
\usepackage{amssymb}
\usepackage{euscript}
\begin{document}
\setlength{\arraycolsep}{.136889em}
\renewcommand{\theequation}{\thesection.\arabic{equation}}
\newtheorem{thm}{Theorem}[section]
\newtheorem{propo}{Proposition}[section]
\newtheorem{lemma}{Lemma}[section]
\newtheorem{corollary}{Corollary}[section]
\newtheorem{remark}{Remark}[section]
\def\begg{\begin{equation}}
\def\endd{\end{equation}}
\def\ep{\varepsilon}
\def\noo{n\to\infty}
\def\al{\alpha}
\def\be{\bf E}
\def\bp{\bf P}
\medskip
\centerline{\Large\bf On the local time of the Half-Plane Half-Comb walk}

\bigskip\bigskip
\bigskip\bigskip
\renewcommand{\thefootnote}{1}
\noindent {\textbf{Endre Cs\'{a}ki}
\newline
Alfr\'ed R\'enyi Institute of Mathematics, 
Budapest, P.O.B. 127, H-1364, Hungary. E-mail address:
csaki.endre@renyi.hu

\bigskip

\noindent {\textbf{Ant\'{o}nia F\"{o}ldes}
\newline
Department of Mathematics, College of Staten Island, CUNY, 2800
Victory Blvd., Staten Island, New York 10314, U.S.A.  E-mail
address: Antonia.Foldes@csi.cuny.edu

\medskip
\noindent{\bf Abstract}
The Half-Plane Half-Comb walk is  a random walk on the plane, when we have
a square lattice on the upper half-plane and a comb structure on the lower
half-plane, i.e., horizontal lines below the $x$-axis are removed. We prove that the probability 
that this walk return to origin in $2N$ steps is asymptotically equal to $2/ (\pi N).$
As a consequence we prove strong laws  and a limit distribution for the local time.  


 \vspace{.1cm}

\section{Introduction and main results}
\renewcommand{\thesection}{\arabic{section}} \setcounter{equation}{0}
\setcounter{thm}{0} \setcounter{lemma}{0}

The properties of a simple symmetric random walk on the square lattice
${\mathbb Z}^2$ have been extensively investigated in the literature
since Dvoretzky and Erd\H os \cite{DE}, and Erd\H os and Taylor \cite{ET}. 
For these and further results we refer to R\'ev\'esz \cite{REV}.

Subsequent investigations concern random walks on other structures of the
plane. For example, a simple random walk on the 2-dimensional
comb lattice that is obtained from ${\mathbb Z}^2$ by removing all
horizontal lines off the $x$-axis was studied by Weiss and Havlin 
\cite{WH}, Bertacchi and Zucca \cite{BERZU}, Bertacchi \cite{BER}, Cs\'aki 
et al. \cite{CSCSFR09}, \cite{CSCSFR11}.

The latter  are particular cases of the so-called anisotropic random walk on the
plane. The general case is given by the transition probabilities
$$ {\bf P}({\bf C}(N+1)=(k+1,j)|{\bf C}(N)=(k,j))=
{\bf P}({\bf C}(N+1)=(k-1,j)|{\bf C}(N)=(k,j))
=\frac{1}{2}-p_j,$$
$$ {\bf P} ({\bf C}(N+1)=(k,j+1)|{\bf C}(N)=(k,j))={\bf P}
({\bf C}(N+1)=(k,j-1)|{\bf C}(N)=(k,j))=p_j,$$
for $(k,j)\in{\mathbb Z^2}$, $N=0,1,2,\ldots$ with
$0<p_j\leq 1/2$ and  $\min_{j\in\mathbb Z} p_j<1/2$. See Seshadri et al.
\cite{SLS}, Silver et al. \cite{SSL}, Heyde \cite{HE} and Heyde et al. 
\cite{HWW}. The simple symmetric random walk corresponds to the case 
$p_j=1/4$, $j=0,\pm1, \pm2,\ldots$, while $p_0=1/4$, $p_j=1/2$, $j=\pm1, 
\pm2,\ldots$ defines  the random walk on the comb.

In  our paper \cite{CSCSFR12}  we combined the simple symmetric random walk 
with a random walk on a comb, when $p_j=1/4$, $j=0,1,2,\ldots$ and
$p_j=1/2$, $j=-1,-2,\ldots$, i.e., we have a square lattice on the upper
half-plane, and a comb structure on the lower half-plane. We call this model
Half-Plane Half-Comb (HPHC) and denote the random walk on it by ${\bf 
C}(N)=(C_1(N),C_2(N)), \, N=0,1,2,\ldots$ 
Here, for  convenient information, we first repeat the precise construction of this walk,    as it was given in [7]:

On a suitable probability space  consider two independent simple symmetric
(one-dimensional) random walks $S_1(\cdot)$, and $S_2(\cdot)$. We may
assume that on the same probability space we have a sequence of
independent geometric random variables $\{Y_i,\, i=1,2,\ldots\}$,
independent from $S_1(\cdot),S_2(\cdot)$, with distribution
\begg
\mathbf{P}(Y_i=k)=\frac{1}{2^{k+1}},\, \, k=0,1,2,\ldots   \label{firstgeo}
\endd
Now horizontal steps will be taken consecutively according to 
$S_1(\cdot)$, and vertical steps consecutively according to $S_2(\cdot)$ 
in the following way. Start from $(0,0)$, take $Y_1$ horizontal steps 
(possibly $Y_1=0$) according to $S_1(\cdot)$, then take 1 vertical step. 
If this arrives to the upper half-plane ($S_2(1)=1$), then take $Y_2$ 
horizontal steps. If, however, the first vertical step is in the negative 
direction ($S_2(1)=-1$), then continue with another vertical step, 
and so on. In  general, if the random walk is on the upper half-plane, 
$(y\geq 0)$ after a vertical step, then take a random number of horizontal 
steps according to the next (so far) unused $Y_j$, independent from the 
previous steps. On the other hand, if the random walk is on the lower half- 
plane ($y<0$) then continue with vertical steps according to $S_2(\cdot)$ 
until it reaches the $x$-axis, and so on.

In  paper  \cite{CSCSFR12} we investigated the almost sure limit properties of this walk 
by using strong approximation methods. 
Our first result was a strong approximation of both components of
the random walk ${\bf C}(\cdot)$ by certain time-changed Wiener processes
(Brownian motions) with rates of convergence. Before stating it, we need
some definitions. Assume that we have two independent standard Wiener
processes $W_1(t),W_2(t),\, \, t\geq 0$, and consider
$$
\alpha_2(t):=\int_0^t I\{W_2(s)\geq 0\}\, ds,
$$
i.e., the time spent by $W_2(\cdot) $ on the non-negative side during the interval
$[0,t]$. The process $\gamma_2(t):=\alpha_2(t)+t$ is strictly increasing,
hence we can define its inverse: $\beta_2(t):=\gamma^{-1}_2(t)$.
Observe that the processes $\alpha_2(t),\, \beta_2(t)$ and $ \gamma_2(t)$
are defined in terms of $W_2(t),$ so they are independent from $W_1(t).$
It can be seen moreover that $0\leq \alpha_2(t)\leq t$, and 
$t/2\leq \beta_2(t)\leq t$. 

\smallskip

{\bf Theorem A } 
{\it On an appropriate probability space for the HPHC  random walk
\newline $\{{\bf C}(N)=(C_1(N),C_2(N));\, N=0,1,2,\ldots\}$
one can construct two independent standard Wiener
processes $\{W_1(t);\, t\geq 0\}$, $\{W_2(t);\, t\geq 0\}$ such that,
as $N\to\infty$, we have with any $\varepsilon>0$ }
$$
|C_1(N)-W_1(N-\beta_2(N))|+|C_2(N)-W_2((\beta_2(N))|
=O(N^{3/8+\varepsilon})\quad {a.s.}
$$
Our second result in  paper \cite{CSCSFR12}   was the following LIL.

{\bf Theorem B }  We have

$$
\limsup_{N\to\infty}\frac{C_1(N)}{\sqrt{N\log\log N}}=\limsup_{N\to\infty}\frac{C_2(N)}{\sqrt{N\log\log N}}=1
\quad a.s.
$$
Furthermore
$$
\liminf_{N\to\infty}\frac{C_1(N)}{\sqrt{N\log\log N}}=-1
\quad a.s., \qquad
\liminf_{N\to\infty}\frac{C_2(N)}{\sqrt{N\log\log N}}=-\sqrt{2}
\quad a.s.
$$

 Moreover  we gave an  explicit formula for the $N$-step  
return probability of the walk, which however was too complicated to conclude the 
asymptotic limit. The aim of the present paper is to study the local time 
of this walk. Based on the just mentioned  formula and a beautiful result of Sparre Andersen,
we first get the asymptotic limit of this return probability, and then use it 
for getting local time results.

\section{Preliminaries}
\renewcommand{\thesection}{\arabic{section}} \setcounter{equation}{0}
\setcounter{thm}{0} \setcounter{lemma}{0}

Let $X_1, X_2...$ be i.i.d. random variables with ${\bf P}(X_1=\pm1)=1/2,$ and  
define $S(0)=0,\,S(i)=\sum_{j=1}^i X_j.$ Then $\{S(n), n=0,1\ldots\}$ is  a 
simple symmetric random walk on the line with local time 
$$\xi(x,n)=\#\{j:0\leq j\leq n, \,\, S(j)=x \},\qquad x\in \mathbb Z,$$ and 
put 
$$A(n)=\sum_{j=0}^\infty\xi(j,n-1),\quad n=1,2,\ldots$$
\begg
{\bf P}(2n,r)={\bf P}(A(2n)=r,S(2n)=0), \quad r=1,2 ,...,2n.  \label{pen}
\endd
Define 
$$ G_n=\# \{j: \, 0\leq j< n, S(j) \geq 0 \} .$$
Then we can rephrase the definition of ${\bf P}(2n,r)$ as follows:
\begg
{\bf P}(2n,r)={\bf P}(G_{2n}=r, S(2n)=0). \label{haha}
\endd
We proved in  \cite{CSCSFR12}, that
\begin{eqnarray}
{\bf P}({\bf C}(2N)&=&(0,0))   \nonumber \\
&=&{2N\choose N}\frac{1}{4^{2N}}+
\sum_{n=1}^N\sum_{r=1}^{2n}{\bf P}(2n,r){2N-2n\choose N-n}\frac{1}{2^{2N-2n}}
{2N-2n+r\choose r}\frac{1}{2^{2N-2n+r}}, \label{main}
\end{eqnarray}
where it was shown that 
\begg {\bf P}(2n,2r-1)={\bf P}(2n,2r),   \label{parity}\endd
and we concluded the following complicated  formula for  ${\bf P}(2n,2r)$ (see Lemma 5.2 
in \cite{CSCSFR12})

$${\bf P}(2n,2r)=\frac{1}{2^{2n}}\sum_{j=1}^{r} \frac{1}{2j-1}{2j-1 \choose j} 
\frac{1}{2n+1-2j}{2n+1-2j \choose n+1-j}.$$
\noindent
However, in order to proceed, we need a closed form for ${\bf P}(2n,2r).$

Sparre Andersen \cite{AES} proved some elegant results about the fluctuation 
of the sums of random variables. We only quote the case of simple symmetric 
random walk  of his much more general results. In his formula (5.12) he defines 
$$ K_n=\# \{j:  0< j  \leq n ,S(j) > 0 \} ,$$
    and  gives the probability of 
$${\bf P}(K_{2n-1}=2r,\,S(2n)=0)={\bf P }(K_{2n-1}=2r+1,\,S(2n)=0)$$
$$=\frac{1}{2}c_{2n} \frac{1}{n+1}\left(1+ \frac{n-2r}{n} 
{-\frac{1}{2}\choose r} 
{-\frac{1}{2}\choose n-r}{-\frac{1}{2}\choose n}^{-1} \right), \quad \quad 
r=0,...,n-1 ,
$$
where $$c_{2n}:={\bf P }(S(2n)=0)=(-1)^n {-\frac{1}{2}\choose n}=
\frac{1}{2^{2n}}{2n\choose n}.$$
  
\begg {\bf P}(K_{2n-1}=2r,\,S(2n)=0)=\frac{1}{2^{2n+1}}{2n\choose n}\frac{1}{n+1}
\left(1+\frac{n-2r}{n}\frac{{2r\choose r}{2n-2r\choose n-r}}
{{2n\choose n}}\right). 
\endd
In the above formulas  for any real number $\alpha$ we used the notation  $\displaystyle{{\alpha \choose r}=\frac{\alpha (\alpha-1)(\alpha-2)\ldots(\alpha-r+1)}{r!}. }$ 

However ${\bf P }(2n,2r)={\bf P }(G_{2n}=2r,\,S(2n)=0),$ given in (\ref{haha}) is 
slightly different from the above one. We will show the following 
\begin{lemma}
\begg {\bf P }(2n,2r)={\bf P }(G_{2n}=2r,\,S(2n)=0)=\frac{1}{2^{2n+1}}{2n\choose n}
\frac{1}{n+1}\left(1+\frac{2r-n}{n}\frac{{2r\choose r}{2n-2r\choose n-r}}
{{2n\choose n}}\right). \label {ander}
\endd
\end{lemma}
\noindent
{\bf Proof}: Recall the definition of $K_n$ and $G_n$ and let

$$M_{n}=\# \{j:  0< j  \leq n, S(j) \leq 0 \}.$$
Then observe that
$${\bf P }(M_{2n}=k, S(2n)=0)={\bf P }(G_{2n}=k, S(2n)=0).$$
Moreover, the following two events are the same:

\begg
\{M_{2n}=r, S(2n)=0\}=\{K_{2n-1}=2n-r, S(2n)=0\}.
\endd

So

\begg 
{\bf P }(2n,2r)={\bf P }(G_{2n}=2r, S(2n)=0)={\bf P }(M_{2n}=2r, S(2n)=0)={\bf P }(K_{2n-1}=2n-2r, 
S(2n)=0),
\endd
which immediately implies our lemma.  $\Box$

\smallskip

Recall now the definition  of  the sequence of i.i.d. geometric random variables given in the introduction

\begg {\bf P}(Y_i=k)=2^{-(k+1)} \quad i=1,2...,\quad k=0,1,2... ,
\endd
and let 
\begg U= U_{K}=\sum_{i=1}^{K} Y_i.   \label{theu}
\endd
 Then  $U_{K}$ is negative binomial with $E(U_{K})=K,$ 
 $Var(U_{K})=2K$ and
 \begg
 {\bf P}(U_{K}=r)={K-1+r\choose r}\frac{1}{2^{K+r}},\,\, r=0,1,2,... 
 \endd
 We will need the following two well-known identities about the negative 
binomial distribution:

\begg
 \sum_{r=0}^{a} {a+r \choose r} \frac{1}{2^{a+r}}=1\endd
 
\begg
 \sum_{r=0}^{\infty} {a+r \choose r} \frac{1}{2^{a+r}}=2\endd
See the first one,  e.g., in Pitman \cite{P} (page 220), while the second 
one is equivalent with 

\noindent
$\sum_{r=0}^{\infty} {\bf P }(U_K=r)=1.$

\bigskip
\noindent
{\bf Lemma A}  {\it Berry-Esseen bound}: \cite {PE} (page 150) 
{\it  Let $X_1, ...X_n$ 
be i.i.d.  random variables. Let
$$ E(X_1)=0,\quad Var(X_1)=\sigma^2>0,\quad  E(|X|^3)<\infty, \quad 
\rho=E(|X|^3)/\sigma^3.$$ 
Then with some constant $A>0$ we have

\begg \sup_x \left|{\bf P}\left(\sigma^{-1}n^{-1/2}\sum_{j=1}^nX_j<x \right) 
-\Phi (x)\right|\quad\leq A\rho\,  n^{-1/2},  \label{ber}\endd
where $\Phi(\cdot)$ is the standard normal distribution function.}
 
\bigskip
\noindent

In what follows we will use a result of Chen \cite{CX}, about Harris recurrent Markov chains, so here we recall his definition.
Let $\{X_n\}_{n\geq 0}$ be a recurrent Markov chain with state space $(E, \mathcal{E}), $  transition probability  $P(x,A)$  and invariant  measure $\mu.$ Recall that
$\{X_{n}\}_{n\geq 0}$ is called Harris  recurrent  if it is irreducible and for any $A \in  \mathcal{E}^+$, and  initial  distribution  $\nu$,
$$P_{\nu}(X_n \in A \quad  {\rm infinitely\, often})=1,$$
where $\mathcal{E}^+=\{A\in \mathcal{E}; \mu(a)>0 \}.$  By Harris recurrence  the invariant  measure $\mu$ uniquely (up to  a constant  multiplier) exists. Obviously our HPHC  walk is Harris recurrent.

\section{Asymptotic return probability}
\renewcommand{\thesection}{\arabic{section}} \setcounter{equation}{0}
\setcounter{thm}{0} \setcounter{lemma}{0}

We want to determine the asymptotic probability that the HPHC  random walk 
returns to the starting point in $2N$ steps. 
\begin{thm} For the asymptotic return probability of the HPHC walk, starting at $(0,0),$  we have
$$
{\bf P }({\bf C}(2N)=(0,0))\sim \frac{2}{\pi N}, \quad as\quad N\to\infty.
$$
\end{thm}

\noindent{\bf Proof:} Recall the definition of $U_K$ in  (\ref{theu}). In 
what follows let $U=U_{2N-2n+1}.$ Introduce the notation
$$ 
Q(r,n):=\frac{{2r\choose r}{2n-2r\choose n-r}}
{{2n\choose n}}.$$
Combining formulas (\ref {main}), (\ref{parity})  and (\ref {ander}) we 
have that

\begin{eqnarray}
&& {\bf P }({\bf C}(2N)=(0,0))={2N\choose N}\frac{1}{4^{2N}}   \label{big} \\
&&+\sum_{n=1}^N {2N-2n\choose N-n} \frac{1}{2^{2N-2n}}\sum_{r=1}^{2n} {\bf P }(2n,r) 
{2N-2n+r  \choose r}\frac{2}{2^{2N-2n+r+1}}  \nonumber \\
&&={2N\choose N}\frac{1}{4^{2N}} + 
\sum_{n=1}^N {2N-2n\choose N-n} \frac{1}{2^{2N-2n}}\sum_{r=1}^{n} 2 {\bf P }(2n,2r) 
{\bf P }(2r-1\leq U \leq 2r) \nonumber \\
&&={2N\choose N}\frac{1}{4^{2N}} +
\sum_{n=1}^N {2N-2n\choose N-n}\frac{1}{2^{2N}} \frac{1}{n+1} 
{2n \choose n}\sum_{r=1}^{n}
\left(1+\frac{2r-n}{n}Q(r,n) \right){\bf P }(2r-1\leq U \leq 2r). \nonumber
\end{eqnarray}
\noindent
The first term in (\ref{big}) is negligible, since 
$$
{2N\choose N}\frac{1}{4^{2N}}=O\left(\frac{1}{4^N}\right).
$$

\noindent
Thus 

\begin{eqnarray} &&{\bf P }({\bf C}(2N)=(0,0)) \nonumber\\
&& \sim \sum_{n=1}^N {2N-2n\choose N-n}\frac{1}{2^{2N}}{2n\choose n} 
\frac{1}{n} {\bf P }(U \leq 2n)\nonumber\\
&&+\sum_{n=1}^N {2N-2n\choose N-n}\frac{1}{2^{2N}}{2n\choose n} 
\frac{1}{n} 
\sum_{r=1}^{n} \frac {2r-n}{n}Q(r,n){\bf P }(2r-1\leq U \leq 2r)= I+II.\nonumber
\end{eqnarray}

Observe that 
\begg{2N-2n\choose N-n}\frac{1}{2^{2N}}{2n\choose n} =c_{2N-2n}c_{2n}\sim 
\frac{1}{\pi}\frac{1}{\sqrt{n}}\,\frac{1}{\sqrt{N-n}}, \quad {\rm when} \quad   
n\to \infty \quad {\rm and}\quad N-n  \to \infty. \label{firstbin}
\endd
Moreover,  if only $n \to \infty$ but  $N-n$ might be  small, then  
\begg c_{2N-2n}\,c_{2n}\leq \frac{c}{\sqrt{n}}. \label{secondbin}
\endd
Here and in what follows $c$ is a positive constant whose value  can 
change from line to line. 
It is  clear that for $1\leq r \leq n$ 
$$-1\leq \frac{2r-n}{n} \leq 1.$$

We will show that term  II  is negligible compared to term I, so we use the 
above fact to give the following upper bound for II:

$$
| \, II \,|\leq II^*:=\sum_{n=1}^N\frac{1}{n}{2n\choose n}
{2N-2n\choose N-n}\frac{1}{2^{2N}}
\sum_{r=1}^nQ(r,n) {\bf P }(2r-1\leq U\leq 2r).
$$

First we deal with the term I, dividing the sum for $n$ into 5 parts:
\begin{eqnarray*}
&(i)&\quad 1\leq n<\frac{N}{4}\\
&(ii)&\quad \frac{N}{4}\leq n< \frac{N}{2}-N^{1/2+\alpha}\\
&(iii)&\quad \frac{N}{2}-N^{1/2+\alpha}\leq n< \frac{N}{2}+N^{1/2+\alpha} \\
&(iv)&\quad \frac{N}{2}+N^{1/2+\alpha}\leq n< N-N^{1/2-\alpha} \\
&(v)& \quad N-N^{1/2-\alpha} \leq n \leq N,
\end{eqnarray*}
with some $0<\alpha<1/2.$ Observe  that
\begg
I=\sum_{n=1}^N \frac{1}{n } c_{2N-2n} \,c_{2n} {\bf P }(U\leq 2n).
\endd
Let us start with $(i).$ In this case we can use the estimation
\begg
{\bf P }(U\leq 2n)=\sum_{r=0}^{2n} {2N-2n+r\choose r}\frac{1}{2^{2N-2n+r+1}}
\leq 2n{2N\choose 2n}\frac{1}{2^{2N}}, \label{first}
\endd
since the largest term in the previous sum corresponds to $r=2n$. Thus 
\begin{eqnarray}
\sum_{(i)}&\leq& \sum_{1\leq n<N/4}  c_{2N-2n}\,c_{2n} \,
\frac{1}{n} 2n{2N\choose 2n}\frac{1}{2^{2N}} \nonumber \\
\leq &c&  \frac{1}{2^{2N}}\sum_{1\leq n<N/4} {2N\choose 2n}\leq 
c \frac{1}{2^{2N}}\frac{N}{4}  {2N \choose N/2} \leq c\sqrt{N} 
\left(\frac{4}{3\sqrt{3}}\right)^N \label{second}
\end{eqnarray}
with some constant $c$, by observing that the first two factor in our sum is 
the product of two probabilities. We used Stirling formula to get the last 
inequality. 

In  case $(ii)$ we use normal approximation for 
negative binomial distribution, with Berry-Esseen bound as in (\ref{ber}) to get 
that for $n$ belonging to the set $(ii)$
$$
{\bf P }(U \leq 2n)=\Phi\left(\frac{4n-2N-1}{\sqrt{2(2N-2n+1)}} \right) 
+O\left(\frac{1}{\sqrt{N-n}}\right)\leq  \Phi(-2N^{\alpha})
+\frac{c}{\sqrt{N}}\leq \frac{c}{\sqrt{N}},
$$
being the normal  term exponentially small. Moreover, using  (\ref{firstbin} )
$$
\sum_{(ii)}\leq c\sum_{N/4< n\leq N/2-N^{1/2+\alpha}}
\frac{1}{  n\sqrt{ n(N-n)}} \frac{1}{\sqrt{N}}
\leq \frac{c}{N^{3/2}}.
$$

Considering now term $(iii),$ we can overestimate  ${\bf P }(U \leq 2n)$ by 1, and 
 obtain, using  (\ref{firstbin}) again,
$$
\sum_{(iii)}\sim \frac{1}{\pi}\sum_{N/2-N^{1/2+\alpha}\leq n
< N/2+N^{1/2+\alpha}}\frac{1}{n^{3/2}(N-n)^{1/2}}\leq \frac{c}{N^{3/2-\alpha}}.
$$

Skipping term $(iv)$ to finish estimating the negligible terms, it is easy to 
see that
$$
\sum_{(v)} \leq \sum_{N-N^{1/2-\alpha}\leq n< N}
 \frac{1}{n^{3/2 } }\leq c \frac{N^{1/2-\alpha}}{N^{3/2}}
=\frac{c}{N^{1+\alpha}}.
$$
using again only that ${\bf P }(U \leq 2n)\leq 1$ and  (\ref{secondbin}).

Now we want to show that part $(iv)$ in  sum I will give the order of 
magnitude claimed in the theorem. It is easy to see  by  normal 
approximation again that for $n \in (iv)$ we obtain 
$$
\Phi(cN^\alpha)\leq {\bf P }(U\leq 2n)\leq 1,
$$
 to conclude that  for $n \in (iv)$ 
 $$
{\bf P }(U\leq 2n)=1-o(1), \quad as \quad N\to\infty.
$$
   So we need the asymptotic value of  
\begg
\sum_{(iv)}\sim \frac{1}{\pi}\sum_{N/2+N^{1/2+\alpha}\leq n
< N-N^{1/2-\alpha}}\frac{1}{n^{3/2}(N-n)^{1/2}}. \label{tricky}
\endd

\noindent
By showing  that   
$$
\sum_{N/2\leq n\leq N/2+N^{1/2+\alpha}}  \frac{1}{n^{3/2}(N-n)^{1/2}}
\leq c\frac{N^{1/2+\alpha}}{N^2}=\frac{c}{N^{3/2-\alpha}}
$$
and 
$$
\sum_{N-N^{1/2-\alpha}\leq n\leq N} \frac{1}{n^{3/2}(N-n)^{1/2}}
\leq c\frac{N^{1/2-\alpha}}{N^{3/2}}=\frac{c}{N^{1+\alpha}},
$$
we can extend the interval of summation in  (\ref{tricky}) without changing 
the limit of the sum as follows:

$$
I \sim \frac{1}{\pi N} \sum_{N/2<n<N}
\frac{1}{\left(\frac{n}{N}\right)^{3/2}\left(1-\frac{n}{N}\right)^{1/2}}
\frac{1}{N}
\sim \frac{1}{\pi N}\int_{1/2}^1\frac{dv}{v^{3/2}(1-v)^{1/2}}
= \frac{2}{\pi N}.
$$

Concerning the term $II^*$, it is clear that  $Q(r,n)$ being a probability, 
the four negligible terms  which we investigated as terms of $I$ are also 
negligible compared to the main term. The only problem is to estimate the 
sum $II^*$  for $n\in$  $(iv)$. This however  is a delicate calculation. We 
split the sum for $r$ into 3 parts:

\ \ \ (1)\ \  $0\leq r\leq n/4$,

\ \ \ (2)\ \  $n/4 <r\leq n-n^{\beta}$,

\ \ \ (3)\ \  $n-n^{\beta}<r\leq n$.

with some $0<\beta<1/2+\alpha<1.$

\smallskip

For (1) we use that $Q(r,n)$ is a probability, obtaining just as in $(i)$ in (\ref{first}) that
$$
\sum_{r \leq n/4}  \,Q(r,n){\bf P }(2r-1\leq U \leq 2r) \leq c \,{\bf P }(U\leq N/2)<\frac{N}{2}{2N\choose N/2}\frac{1}{2^{2N}},
$$ 
So
\begin{eqnarray*}
&&\sum_{n\in (iv)} \frac{1}{n}{2n\choose n}{2N-2n\choose N-n}\frac{1}{2^{2N}}
\sum_{r\in (1)}Q(r,n) {\bf P }(2r-1< U\leq 2r) \\
&& 
\leq \sum_{ n\in(iv)} \frac{cN}{n}    
{2N\choose N/2}\frac{1}{2^{2N}} \leq c    {2N\choose N/2}\frac{1}{2^{2N}} \leq c \left(\frac{4}{3\sqrt{3}}\right)^N,
\end{eqnarray*} where the last inequality  is coming from Stirling  formula as in (\ref{second}).

In case (2), using Stirling formula, we have 
\begin{eqnarray}
&&\sum_{r\in (2)}Q(r,n) {\bf P }(2r-1\leq U \leq 2r)
\leq  \sum_{r\in (2)}\frac{c\sqrt{n}}{\sqrt{r}\sqrt{n-r}}
{\bf P }(2r-1\leq U \leq 2r) \nonumber\\
&&\leq\frac{c}{n^{\beta/2}}\sum_{r\in (2)}{\bf P }(2r-1\leq U \leq 2r)\leq \frac{c}
{n^{\beta/2}}\leq\frac{c}{N^{\beta/2}},\nonumber
\end{eqnarray}
where the last inequality holds as  $n\in(iv).$
Consequently, similarly to $(iv)$ in calculating I, we have 
$\frac{c}{N^{\beta/2}}$ times the sum in (\ref{tricky}) implying that
$$
\sum_{n\in (iv)}\sum_{r\in (2)}\leq \frac{c}{N^{1+\beta/2}}.
$$

For the case $r \in$ (3) we have

$$
\sum_{r\in (3)}Q(r,n) {\bf P }(2r-1\leq U \leq 2r)\leq 
c{\bf P }(2n-2n^{\beta}-1\leq U\leq 2n)\leq c {\bf P }(U\geq 2n-2n^{\beta}-1).
$$

Recall now that $U=U_{2N-2n+1}$ with  $E(U)=2N-2n+1$ and $Var(U)=2(2N-2n+1).$

Applying now Chebyshev inequality in the form

$${\bf P}(X-\mu\geq x\sigma)\leq {\bf P}(|X-\mu|\geq x\sigma)\leq\frac{1}{x^2}$$

we arrive to 
\begin{eqnarray*}
&&{\bf P}(U-2N+2n-1\geq 4n-2N-2n^\beta-2)=
{\bf P}\left( U-E(U)\geq \frac{4n-2N-2n^\beta-2}{(4N-4n+2)^{1/2}} \sigma \right ) \nonumber \\
&&\leq \frac{4N-4n+2}{(4n-2N-2n^\beta-2)^2 }\sim \frac{N-n}{(2n-N-n^{\beta})^2} .\nonumber \\
\end {eqnarray*}
Being  $n\in(iv)$ we have $\frac{N}{2}+N^{1/2+\alpha}\leq n< N-N^{1/2-\alpha}$, implying that
$$N-n\leq N/2 \quad {\rm and}  \quad 2n-N\geq 2 N^{1/2+\alpha}.$$
Knowing  also that $1>1/2+\alpha>\beta>0$  we can conclude that 
$${\bf P }(U\geq 2n-2n^{\beta}-1)\leq \frac{c}{N^{2\alpha}}$$
which goes to zero as $N\to +\infty,$
so the term $\sum_{n\in (iv)}\sum_{r\in (3)}$ is negligible 
compared to  $\sum_{n\in (iv)}$ in the main term.

This completes the proof of  Theorem 3.1. $\Box$

\section{Laws of the iterated logarithm for the local time}
\renewcommand{\thesection}{\arabic{section}} \setcounter{equation}{0}
\setcounter{thm}{0} \setcounter{lemma}{0}

Define the local time of the random walk on the  HPHC  lattice as

$$ \Xi((k,j) ,N)=\sum_{r=0}^N I\{{\bf C}(r)=(k,j)\}, \,\,(k,j) \in 
{\mathbb Z}^2$$
From Theorem 3.1 we can calculate the truncated Green function  $g(\cdot)$ :

$$g(N)=\sum_{k=0}^{[N/2]} {\bf P} ({\bf C}(2k)=0) \sim \frac{2}{\pi} \log N \qquad as\,\, N\to\infty.$$

 Our random walk  being Harris recurrent, we can infer (e.g. Chen \cite{CX}) that 

$$
\lim_{N\to\infty}\frac{\Xi((k_1,j_1),N)}{\Xi((k_2,j_2),N)}=
\frac{\mu(k_1,j_1)}{\mu(k_2,j_2)}\quad {a.s.},
$$
where $\mu(\cdot)$ is an invariant measure. Here the invariant measure is 
defined as the solution of the equation
$$
\mu(A)=\sum_{(k,j)\in{\mathbb Z}^2}\mu(k,j)\mathbf{P}(\mathbf{C}(N+1)\in
A|\mathbf{C}(N)=(k,j)).
$$

For $(k,j)\in {\mathbb Z}^2$, in our case we have
$$
\mu(k,j)=\mu(k+1,j)\left(\frac12-p_j\right)
+\mu(k-1,j)\left(\frac12-p_j\right)
+\mu(k,j+1)p_{j+1}+\mu(k,j-1)p_{j-1},$$

\noindent
where 
$$p_j=\frac{1}{4} \quad {\rm if} \quad  j\geq 0 \quad {\rm and} 
\quad p_j=\frac{1}{2} \quad  {\rm if} \quad j< 0.$$
It is easy to see that
$$
\mu(k,j)=\frac{1}{p_j}, \quad (k,j) \in {\mathbb Z}^2
$$
is one possible invariant measure.
So from Theorem 17.3.2 of Meyn and Tweedie \cite{MT} we get the following  
result
\begin{corollary} For all  integers $k_1, k_2$ we have
$$
 \lim_{N\to\infty}\frac{\Xi((k_1,j_1),N)}{\Xi((k_2,j_2),N)}=
1 \quad {a.s.} \quad {\rm if} \quad j_1\geq 0,\,j_2\geq 0\quad {\rm or}  
\quad j_1< 0,\,j_2< 0 , 
$$
and 
$$
\lim_{N\to\infty}\frac{\Xi((k_1,j_1),N)}{\Xi((k_2,j_2),N)}=2  
\quad {a.s.} \quad  {\rm if} \quad j_1\geq 0,\,j_2< 0.
$$
\end{corollary}
Using $g(N)$ again, we get from Darling and Kac \cite{DK} the following result.
\begin{corollary}

$$\lim_{N\to \infty}{\bf P}\left (\frac{\Xi((0,0),N)}{g(N)} 
\geq x \right)=\lim_{N\to \infty}{\bf P}\left 
(\frac{\pi\,\Xi((0,0),N)}{2\,\log N} \geq x \right)=e^{-x}.$$
\end{corollary}

As to the law of the  iterated logarithm, we get from Theorem 2.4 of Chen 
\cite{CX} that it reads as follows.

\begin{corollary}
$$ \limsup_{N\to \infty} \frac{\Xi((0,0),N)}{\log N \log\log\log N} 
=\frac{2}{\pi} \qquad a.s. $$
\end{corollary}

\smallskip

 To conclude we would like to discuss how   these  results relate to the corresponding  ones for other anisotropic planar walks. In  the anisotropic  walk in general,  everything is determined by the return probability to zero, which allow us to  calculate the Green function, which leads to the results about the local time. As we will see  below the return probability, which we got for the HPHC walk is only differ in a constant from the return probability of the simple symmetric walk of the plane, and much smaller than the return probability to   zero of the two dimensional comb. 
 
 As far as we know the return probability  to zero for the anisotropic random walk is not known. However for the periodic  anisotropic random walk  ${\bf C^P}(\cdot) $ which is defined by
$p_j=p_{j+L} $ for each $j\in {\mathbb Z},$ where L is a positive integer, we proved  in \cite{CSCSFR13} that 

 $${\bf P}({\bf C^P}(2N) =(0,0))\sim\frac{1}{4\pi Np_0\sqrt{\gamma-1}} \quad as \quad   N\to \infty,$$
where
$$ \gamma=\frac{\sum_{j=0}^{L-1}p^{-1}_j}{2L}.$$  
This leads to the following local time results

$$
 \lim_{N\to\infty}\frac{\Xi^P((0,0),N)}{\Xi^P((k,j),N)}=\frac{p_j}{p_0}  \qquad a.s. 
$$
$$ \limsup_{N\to \infty} \frac{\Xi^P((0,0),N)}{\log N \log\log\log N} 
=\frac{1}{4p_0\pi\sqrt{\gamma-1}} \qquad a.s. $$

In  case of the simple symmetric walk on the plane is,  when $p_j=p_0=1/4$  for $  j=\pm1, \pm2,\ldots$ and $L=1$
this reduces to  the well known   asymptotic formula
 $${\bf P}({\bf C}(2N) =(0,0))\sim\frac{1}{\pi N} \quad as \quad   N\to \infty.$$
 This leads to the famous Erd\H os -Taylor integral  test \cite{ET}  (see e.g in  \cite{REV}) containing  e.g. that
 
 $$ \limsup_{N\to \infty} \frac{\Xi((0,0),N)}{\log N \log\log\log N} 
=\frac{1}{\pi} \qquad a.s. $$
On the other hand  in the case of the 2-dimensional comb, when $p_0=1/4,$ and $p_j=1/2 $ for $ j=\pm1, \pm2, \ldots$ 
we have from Bertacchi and Zucca \cite{BERZU99} that
  $${\bf P}({\bf C}(2N) =(0,0))\sim\frac{1}{2^{9/2}\Gamma(1/4) N^{3/4}} \quad as \quad   N\to \infty.$$
    For the local time of the two-dimensional comb  \cite{CSCSFR11}, we have  for any fix $k$ 
  
  $$ \limsup_{N\to \infty} \frac {\Xi((k,0),N)}{ N^{1/4} (\log\log N)^{3/4}} =\frac{2^{9/4}}{3^{3/4}}   \qquad a.s. $$
 and for any fix $k$ and any fix $ j \neq 0$ 
   $$ \limsup_{N\to \infty} \frac {\Xi((k,j),N)}{ N^{1/4} (\log\log N)^{3/4}} =\frac{2^{5/4}}{3^{3/4}}   \qquad a.s. $$

As we mentioned above, the local time behavior is determined by the return probability to zero. The order of magnitude (apart from a constant factor) of the return probability for the simple symmetric random  walk,   for the periodic walk  on  the plane discussed above, and for the HPHC walk are the same. So their local time behavior are the same as well. However as the order of the return probability of the comb is different, its local time behavior is very different from the other  three cases above.  It would be interesting to find examples which shed some light upon this transition between these two types of behavior.

\bigskip
{\bf Acknowledgements} We wish to thank the referee of our submission, for careful reading our manuscript,
and for  making a number of helpful suggestions which certainly improved the presentation of this paper.

\end{document}